\documentclass[letterpaper, 12 pt]{article}

\usepackage{hyperref}
\usepackage{blindtext}
\usepackage{tikz}
\usepackage{calc}
\usepackage[top=1in,bottom=1in,left=1in,right=1in,includehead]{geometry}

\usepackage{amsmath,amsthm,amssymb}
\usepackage{graphicx}
\usepackage{caption}

\usepackage[sc]{mathpazo}
\usepackage[T1]{fontenc}

\usepackage[capitalize,nameinlink,noabbrev]{cleveref}

\usepackage{authblk}

\let\set\mathbb

\def\id{\operatorname{id}}
\def\lclm{\operatorname{lclm}}

\title{Quadrant Walks Starting Outside the Quadrant}

\date{\vspace{-6ex}}

\author{Manfred Buchacher, Manuel Kauers, and Am\'{e}lie Trotignon
\thanks{manfred.buchacher@jku.at, manuel.kauers@jku.at, amelie.trotignon@jku.at. All authors were supported by the Austrian FWF grant F5004}}

\affil{\small Institute for Algebra, Johannes Kepler University Linz, Austria}

\usepackage[backend=bibtex]{biblatex}
\begin{document}

\maketitle

\abstract{
We investigate a functional equation which resembles the functional equation for the generating function of a lattice walk model for the quarter plane.
  The interesting feature of this equation is that its orbit sum is zero while its solution is not algebraic.
  The solution can be interpreted as the generating function of lattice walks in $\mathbb{Z}^2$ starting at $(-1,-1)$ and subject to the restriction that
  the coordinate axes can be crossed only in one direction.
  We also consider certain variants of the equation, all of which seem to have transcendental solutions.
  In one case, the solution is perhaps not even D-finite. 
}


\section{Introduction}

The investigation of lattice walks with small steps restricted to a quadrant has made astonishing progress during the past
years~\cite{bousquet10,fayolle99,bostan14b,kurkova15,bostan16b,courtiel17,bernardi17,bostan17,dreyfus18}.
The central problem in this context is to decide for a given step set $S\subseteq\{-1,0,1\}^2\setminus\{(0,0)\}$ whether the
generating function $F(x,y,t)=\sum_{n=0}^\infty\sum_{i,j=0}^n a_{i,j,n} x^i y^j t^n$ counting the number $a_{i,j,n}$ of walks
in $\set N^2$ starting at $(0,0)$, ending at $(i,j)$, and consisting of exactly $n$ steps, each step taken from~$S$, is D-finite.
If so, it is further of interest whether it is even algebraic. 
Although it is not obvious at first glance, it is meanwhile well understood how the finiteness of a certain group associated to
the model implies the D-finiteness of the generating function, and how the zeroness of the so-called orbit sum
associated to the model implies its algebraicity~\cite{fayolle10,kurkova12,raschel12}. 

For simplicity, let us focus on the step set $S=\{\leftarrow,\rightarrow,\uparrow,\downarrow\}$. If $Q(x,y,t)$ is the generating
function for this model, the combinatorial definition translates into the functional equation
\[
(1 - (x+y+\bar x+\bar y)t) Q(x,y,t) = 1 - \bar x t Q(0,y,t) - \bar y t Q(x,0,t),
\]
where we write $\bar x=1/x$ and $\bar y=1/y$ for short.
The group of the model can be used to solve this equation.
In the present example, it is generated by the rational maps $\Phi=(\binom xy\mapsto\binom{\bar x}y)$ and $\Psi=(\binom xy\mapsto\binom x{\bar y})$
under composition, i.e., $G=\{\id,\Phi,\Psi,\Phi\circ\Psi\}$.
The idea for solving the functional equation is to let the elements of the group act on it to get four copies of the equation, and then
take a linear combination of these four copies with the aim of canceling the terms $Q(\cdots)$ appearing on the right. This leads to
\begin{alignat*}1
  &(1 - (x+y+\bar x+\bar y)t) \bigl(xyQ(x,y,t) - \bar xy Q(\bar x,y,t) - x\bar y Q(x,\bar y,t) + \bar x\bar y Q(\bar x,\bar y,t)\bigr)\\
  &\qquad{}= xy - \bar x y-x\bar y+\bar x\bar y.
\end{alignat*}
The expression on the right is the \emph{orbit sum.}
Divide by $1-(x+y+\bar x+\bar y)t$ and observe next that $xyQ(x,y,t)$ is the only term on the left whose exponents with respect to $x$ and $y$ are positive,
while for all terms $x^iy^jt^n$ appearing in any of the other terms on the left we have $i<0$ or $j<0$. 
Therefore, by extracting the positive part, we can eliminate the unwanted terms $Q(\bar x,y,t),Q(x,\bar y,t),Q(\bar x,\bar y,t)$ and get
\begin{equation}\label{eq:Q}\tag{Q}
  x y Q(x,y,t) = [x^>][y^>] \frac{xy - \bar x y-x\bar y+\bar x\bar y}{1 - (x+y+\bar x+\bar y)t}.
\end{equation}
Since extracting the positive part preserves D-finiteness, it follows that $Q(x,y,t)$ is D-finite.

The step sets
$\{\leftarrow,\downarrow,\nearrow\}$ (Kreweras),
$\{\rightarrow,\uparrow,\swarrow\}$ (reverse Kreweras),
$\{\leftarrow,\rightarrow,\downarrow,\uparrow,\penalty0\nearrow,\penalty0\swarrow\}$ (double Kreweras), and
$\{\leftarrow,\swarrow,\rightarrow,\nearrow\}$ (Gessel)
also have a certain finite group of rational maps associated to them, but the approach above for solving the functional equations for the
generating functions fails, because the orbit sum turns out to be zero in these cases.
Using more sophisticated arguments, it can be shown that the generating functions for these models are not only D-finite but in fact
algebraic~\cite{bousquet10,bostan10}. In fact, the generating function happens to be algebraic if and only if the orbit sum vanishes. 

The equivalence between zeroness of the orbit sum and algebraicity of the generating function is not an accident, but it can be explained~\cite[Sec.~8 and Sect.~9.1]{kurkova15}.
However, as we shall show in this paper, the equivalence does not hold in all circumstances.
Consider the following slight variation of the functional equation quoted above for the step set $\{\leftarrow,\rightarrow,\uparrow,\downarrow\}$:
\begin{equation}\label{eq:F}\tag{F}
(1 - (x+y+\bar x+\bar y)t) F(x,y,t) = \bar x\bar y - \bar x t F(0,y,t) - \bar y t F(x,0,t).
\end{equation}
The only difference is that we replaced the inhomogeneous term $1$ by $\bar x\bar y$. If we now multiply the equation by $xy$, as above,
then let the four group elements act on the equation, as above, and then take the weighted sum of the resulting equations, as above, we obtain
\[
(1 - (x+y+\bar x+\bar y)t) \bigl(xyF(x,y,t) - \bar xy F(\bar x,y,t) - x\bar y F(x,\bar y,t) + \bar x\bar y F(\bar x,\bar y,t)\bigr)=0,
\]
and it is not clear how to proceed from here.

It is clear that any solution $F(x,y,t)\in\set Q[x,\bar x,y,\bar y][[t]]$ will have the form $\bar x\bar yt^0+\cdots$, so before we proceed,
we should clarify what we mean by the expressions $F(0,y,t)$ and $F(x,0,t)$ on the right hand side of~\eqref{eq:F}.
There are several options. 
A natural interpretation is $F(0,y,t)=[x^0]F(x,y,t)$ and $F(x,0,t)=[y^0]F(x,y,t)$. We consider this case in Section~\ref{sec:all}.
Other interpretations also include certain restrictions on the other variable.
For example, we could choose to read $F(0,y,t)$
as $[x^0][y^<]F(x,y,t)$ (keeping only negative exponents of~$y$),
as $[x^0][y^\geq]F(x,y,t)$ (keeping only nonnegative exponents of~$y$),
as $[x^0][y^\leq]F(x,y,t)$ (keeping only nonpositive exponents of~$y$), or
as $[x^0][y^>]F(x,y,t)$ (keeping only positive exponents of~$y$),
and analogously for $F(x,0,t)$. 
For some of these interpretations, we can show that the solution of the functional equation is D-finite.
For some of them, we can show that the correctness of a guessed differential equation for the specialization $F(1,1,t)$ of the solution $F(x,y,t)$
implies the transcendence of the solution (Sections \ref{sec:nnp} and~\ref{sec:np/npp}).
One case seems to be more complicated. In this case, we conjecture that the solution is not even D-finite (Section~\ref{sec:pp}).

\section{Four Compartments}\label{sec:all}

With the interpretation $F(0,y,t)=[x^0]F(x,y,t)$ and $F(x,0,t)=[y^0]F(x,y,t)$, the solution $F(x,y,t)$ of the functional equation~\eqref{eq:F}
counts walks that start at $(-1,-1)$ and move through the plane $\set Z^2$ subject to the restriction that the axes of the coordinate system
can be passed only in one direction (west to east or south to north, respectively).
We claim that the generating function $F(x,y,t)$ counting walks in this model is D-finite. To show this, define 
$F_1=[x^<][y^<]F$, $F_2=[x^<][y^\geq]F$, $F_3=[x^\geq][y^<]F$, and $F_4=[x^\geq][y^\geq]F$, so that
$F=F_1+F_2+F_3+F_4$.
\begin{center}
\begin{tikzpicture}[scale=.33]
  \draw[->](-5.5,0)--(5.5,0);
  \draw[->](0,-5.5)--(0,5.5);
  \foreach \x in {-5,...,5} \foreach \y in {-5,...,5} \fill (\x,\y) circle(2pt);
  \draw[blue] (-1,-1) circle(3pt);
  \begin{scope}[red,thick]
  \foreach \x in {-3,3}
  \foreach \y in {-3,3}
  \foreach \u/\v in {1/0,0/1,-1/0,0/-1}
  \draw[->] (\x,\y) --+(\u,\v);
  \foreach \x in {-3,3}
  \foreach \u/\v in {1/0,0/1,-1/0} {
    \draw[->] (\x,0) --+(\u,\v);
    \draw[->] (0,\x) --+(\v,\u);
  }
  \draw[->] (0,0)--(1,0);
    \draw[->] (0,0)--(0,1);
  \end{scope}
\end{tikzpicture}\hfil
  \begin{tikzpicture}[scale=.33]
    \fill[lightgray,rounded corners=2pt] (-.67,-.67) rectangle (-5.25,-5.25)
    (-.67,-.33) rectangle (-5.25,5.25) (-.33,-.67) rectangle (5.25,-5.25) (-.33,-.33) rectangle (5.25,5.25);
    \fill[lightgray] (-2,-2) rectangle (-5.25,-5.25) (-2,2) rectangle (-5.25,5.25) (2,2) rectangle (5.25,5.25) (2,-2) rectangle (5.25,-5.25);
    \draw[->](-5.5,0)--(5.5,0);
    \draw[->](0,-5.5)--(0,5.5);
    \foreach \x in {-5,...,5} \foreach \y in {-5,...,5} \fill (\x,\y) circle(2pt);
    \draw[blue] (-1,-1) circle(3pt);
    \draw[thick,rounded corners=2pt] 
    (-5.25,-.67)--(-.67,-.67)--(-.67,-5.25) (-2.5,-2.5) node {$F_1$}
    (-5.25,-.33)--(-.67,-.33)--(-.67,5.25) (-2.5,2.5) node {$F_2$}
    (-.33,-5.25)--(-.33,-.67)--(5.25,-.67) (2.5,-2.5) node {$F_3$}
    (-.33,5.25)--(-.33,-.33)--(5.25,-.33) (2.5,2.5) node {$F_4$};
  \end{tikzpicture}
\end{center}
The equation for $F$ translates into the following system of functional equations for $F_1,F_2,F_3,F_4$,
where we write $S=x+y+\bar x+\bar y$:
\begin{alignat*}5
  F_1 &= \bar x\bar y &&+ S t F_1 && - t[\bar x] F_1 - t[\bar y] F_1\\ 
  F_2 &= t [\bar y] F_1 &&+ S t F_2 &&- t[\bar x] F_2 - \bar yt [y^0] F_2\\ 
  F_3 &= t [\bar x] F_1 &&+ S t F_3 &&- t[\bar y] F_3 - \bar xt [x^0] F_3\\ 
  F_4 &= \underbrace{\rule[-5pt]{0pt}{0pt}t [\bar x] F_2 + t [\bar y] F_3}_{\text{\mathstrut initial conditions}}&& +
         \underbrace{\rule[-5pt]{0pt}{0pt}S t F_4}_{\text{\kern-2em\mathstrut recurrence\kern-2em}} &&
  \underbrace{\rule[-5pt]{0pt}{0pt}-\bar xt [x^0] F_4 - \bar yt [y^0] F_4}_{\text{\mathstrut boundary conditions}}. 
\end{alignat*}
The equation for $F_1$ does not depend on $F_2,F_3,F_4$ and can therefore be solved directly.
In fact, we have $F_1(x,y,t)=\bar x\bar y Q(\bar x,\bar y,t)$ for the $Q(x,y,t)$ from equation~\eqref{eq:Q}.

Knowing~$F_1$, we can solve the second equation for~$F_2$ by the same technique.
Noting that $[\bar y]F_1$ is independent of~$y$, the result is
\begin{alignat*}1
  F_2(x,y,t) &= \bar y[x^<][y^>]\frac{t (y [\bar y]F_1(x,y,t) - y [\bar y]F_1(\bar x,y,t) - \bar y[\bar y]F_1(x,y,t) + \bar y [\bar y]F_1(\bar x,y,t))}{1 - S t}\\
  &= t\bar y[x^<][y^>]\frac{(y-\bar y)[\bar y](F_1(x,y,t) - F_1(\bar x,y,t))}{1 - S t},
\end{alignat*}
so $F_2$ is D-finite because it is the positive/negative part of a D-finite series.
Moreover, using $F_1(x,y,t)=\bar x\bar y Q(\bar x,\bar y,t)$ and~\eqref{eq:Q} we get
\[
  F_1(x,y,t) - F_1(\bar x,y,t) = [y^<] \frac{xy-\bar xy-x\bar y+\bar x\bar y}{1- St},
\]
which can be used to simplify the expression for $F_2(x,y,t)$ further to
\[
  F_2(x,y,t) = t\bar y[x^<]\Bigl(\bigl([y^>]\frac{y-\bar y}{1 - St}\bigr)\bigl([\bar y] \frac{xy-\bar xy-x\bar y+\bar x\bar y}{1- St}\bigr)\Bigr).
\]
Because of symmetry, we have $F_3(x,y,t)=F_2(y,x,t)$, so this one is D-finite too, and we can directly proceed to the equation for~$F_4$,
which we can now solve in terms of the known functions $F_2,F_3$, again by letting the group elements act, forming a weighted sum,
dividing by $1-St$ and extracting the positive part. The result is
\[
F_4(x,y,t) = \bar x\bar yt[x^>][y^>] \frac{G(x,y,t)-G(\bar x,y,t)-G(x,\bar y,t)+G(\bar x,\bar y,t)}{1 - St},
\]
with $G(x,y,t)=x y [\bar x] F_2(x,y,t) + x y [\bar y] F_3(x,y,t)$.
We already see at this point that $F_4$ is D-finite, because it is the positive part of a D-finite
series, so we can conclude that $F=F_1+F_2+F_3+F_4$ is D-finite, because it is the sum of four D-finite series.
Moreover, using the expression for $F_2$ derived above, we can state $F_4$ explicitly as
\begin{alignat*}1
F_4(x,y,t) = {}&\bar x\bar y t^2 [y^>]
\Bigl(
\bigl([\bar x]\frac{(y-\bar y)[\bar y]\displaystyle\frac{xy-\bar xy-x\bar y+\bar x\bar y}{1-St}}{1-St}\bigr)
\bigl([x^>]\frac{x-\bar x}{1-St}\bigr)
\Bigr)\\
+{}&\bar x\bar y t^2 [x^>]
\Bigl(
\bigl([\bar y]\frac{(x-\bar x)[\bar x]\displaystyle\frac{xy-\bar xy-x\bar y+\bar x\bar y}{1-St}}{1-St}\bigr)
\bigl([y^>]\frac{y-\bar y}{1-St}\bigr)
\Bigr)
.
\end{alignat*}
The expressions we found for $F_1,F_2,F_3,F_4$ are small enough that we succeeded to use the techniques from~\cite{bostan16b} and Koutschan's package~\cite{koutschan10c}
to construct a certified annihilating operator for~$F(1,1,t)$. We suppress the computational details here and refer the interested reader to the
Mathematica session posted on the accompanying website of this paper~\cite{buchacher20a}. The bottom line is that the coefficient sequence of $F(1,1,t)$ satisfies
the recurrence
\begin{alignat*}1
  &(2 + n)(4 + n)(6 + n)(-1 + 2n + n^2)a_{n+2}\\
  &- 4(3 + n)(-18 + 4n + 9n^2 + 2n^3)a_{n+1}\\
  &-16(1 + n)(2 + n)(3 + n)(2 + 4n + n^2)a_n = 0.
\end{alignat*}
This recurrence has only asymptotic solutions of the form $4^n n^{-1}$ and $(-4)^n n^{-3}$ (as can be found using~\cite{kauers11e,kauers14b}).
Neither of these forms can arise from an algebraic series, so $F(1,1,t)$ must be transcendental.

\bigskip
In summary, we have shown that for the interpretation $F(0,y,t)=[x^0]F(x,y,t)$ and $F(x,0,t)=[y^0]F(x,y,t)$, the solution $F(x,y,t)$ of the functional equation~\eqref{eq:F}
is D-finite but not algebraic.

\section{A Large and a Small Compartment}\label{sec:nnp}

We now turn to the variant of \eqref{eq:F} in which $F(0,y,t)$ is interpreted as $[x^0][y^\geq] F(x,y,t)$, and $F(x,0,t)$ likewise.
In this case, the equation describes a model in which only the nonnegative part of each axis forms a semipermeable barrier for the
walks.
Walks in this model can freely move around in the complement of the north-east quadrant, which is a three quarter plane, and once they
leave this area, they are locked in the north-east quadrant.
It is therefore natural to write the generating function $F(x,y,t)$ for this model as $F=F_1+F_2$
where $F_1=[x^<]F + [x^\geq][y^<]F$ keeps track of the three quarter plane and $F_2=[x^\geq][y^\geq]F$ takes care of the remaining quarter plane.
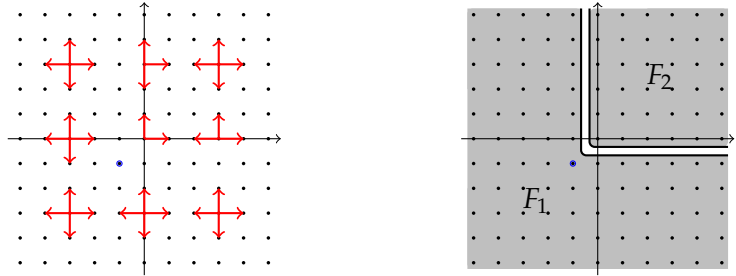
\begin{center}
\begin{tikzpicture}[scale=.33]
  \draw[->](-5.5,0)--(5.5,0);
  \draw[->](0,-5.5)--(0,5.5);
  \foreach \x in {-5,...,5}
  \foreach \y in {-5,...,5}
  \fill (\x,\y) circle(2pt);
  \draw[blue] (-1,-1) circle(3pt);
  \begin{scope}[red,thick]
  \foreach \x in {-3,3}
  \foreach \y in {-3,3}
  \foreach \u/\v in {1/0,0/1,-1/0,0/-1}
  \draw[->] (\x,\y) --+(\u,\v);
  \foreach \x in {-3,3}
  \foreach \u/\v in {1/0,0/1,-1/0} {
    \draw[->] (\x,0) --+(\u,\v);
    \draw[->] (0,\x) --+(\v,\u);
  }
  \draw[->] (-3,0)--+(0,-1);
  \draw[->] (0,-3)--+(-1,0);
  \draw[->] (0,0)--(1,0);
  \draw[->] (0,0)--(0,1);
  \end{scope}
\end{tikzpicture}\hfil
  \begin{tikzpicture}[scale=.33]
    \fill[lightgray,rounded corners=2pt] (-.67,-.67) -- (-.67,5.25) -- (-5.25,5.25) -- (-5.25,-5.25) -- (5.25,-5.25) -- (5.25,-.67) -- cycle
    (-.33,-.33) rectangle (5.25,5.25);
    \fill[lightgray] (-5.25,2) rectangle (-.67,5.25) (-5.25,-2) rectangle (5.25,-5.25) (2,-2) rectangle (5.25,-.67) (2,2) rectangle (5.25,5.25);
    \draw[->](-5.5,0)--(5.5,0);
    \draw[->](0,-5.5)--(0,5.5);
    \foreach \x in {-5,...,5} \foreach \y in {-5,...,5} \fill (\x,\y) circle(2pt);
    \draw[blue] (-1,-1) circle(3pt);
    \draw[thick,rounded corners=2pt]
    (-.67,5.25)--(-.67,-.67)--(5.25,-.67) (-2.5,-2.5) node {$F_1$}
    (-.33,5.25)--(-.33,-.33)--(5.25,-.33) (2.5,2.5) node {$F_2$};
  \end{tikzpicture}
\end{center}
It is known~\cite{bousquet15,raschel19} that the generating function $C(x,y,t)$ for simple walks avoiding the positive quadrant is D-finite.
Hence also $F_1(x,y,t)=\bar x\bar y C(x,y,t)$ is D-finite.
The series $F_2$ counts walks in the quarter plane with initial conditions prescribed by the sections of~$F_1$:
\[
F_2 = \underbrace{\rule[-5pt]{0pt}{0pt}t [\bar x][y^\geq] F_1 + t [\bar y][x^\geq] F_1}_{\text{\mathstrut initial conditions}} +
      \underbrace{\rule[-5pt]{0pt}{0pt} S t F_2}_{\text{\mathstrut\kern-2emrecurrence\kern-2em}}
  \underbrace{\rule[-5pt]{0pt}{0pt}-\bar xt [x^0] F_2 - \bar yt [y^0] F_2}_{\text{\mathstrut boundary conditions}}.
\]
This is again a functional equation which we can solve like in the introduction, the result being a positive part expression in terms of~$F_1$:
\[
F_2 = \bar x\bar y t[x^>][y^>]\frac{
 H(x,y,t) + H(y,x,t)}{1-St},
\]
with $H(x,y,t) = (\bar x-x)[\bar x]\bigl(\bar y[y^\leq]F_1(x,\bar y,t) - y[y^\geq]F_1(x,y,t)\bigr)$.
This implies that $F_2$ is D-finite, so $F=F_1+F_2$ is D-finite as well. 

Unfortunately, in this case we were not able to derive a certified annihilating operator for $F(1,1,t)$ from this expression.
Owing to the size of the equations describing $F_1$, the computations were too expensive.
However, it is easy to guess an annihilating operator $L$ from the first few terms of $F(1,1,t)$.
We found a convincing candidate of order 11 with polynomial coefficients of degree~89.
It is posted on our website~\cite{buchacher20a}.
Assuming that this guessed operator $L$ is correct, we can show that $F(1,1,t)$ is transcendental.
The key observation is that the guessed operator can be written as $L=\lclm(L_1,\dots,L_6)$, where $L_1,\dots,L_6$ are certain irreducible
operators. The factors $L_1,\dots,L_6$ are also posted on our website~\cite{buchacher20a}.
It turns out that $L_1,\dots,L_5$ only have algebraic solutions, while $L_6$ has a logarithmic singularity and therefore can not have any nonzero algebraic solution. The factorization $L=\lclm(L_1,\dots,L_6)$ means that we have $F(1,1,t)=f_1+\cdots+f_6$
where each $f_i$ is a solution of $L_i$.
Since $f_1,\dots,f_5$ are algebraic and $f_6$ is not, it follows that $F(1,1,t)$ is not algebraic unless the term $f_6$ is zero.
In this case however $F(1,1,t)$ would already be annihilated by $\lclm(L_1,\dots,L_5)$, and it can be checked that this
is not the case.

\bigskip
In summary, we have shown that for the interpretation $F(0,y,t)=[x^0][y^\geq] F(x,y,t)$ and $F(x,0,t)=[y^0][x^\geq] F(x,y,t)$,
the solution $F(x,y,t)$ of the functional equation~\eqref{eq:F} is D-finite. Moreover, under the hypothesis that a guessed
annihilating operator for $F(1,1,t)$ is correct, we can also show that $F(1,1,t)$ is transcendental. 

\section{A Small and a Large Compartment}\label{sec:np/npp}


We now turn to the two variants of \eqref{eq:F} where $F(x,0,t)$ and $F(0,y,t)$ are interpreted as $[x^<][y^0] F(x,y,t)$ and $[x^0][y^<] F(x,y,t)$, and as $[x^\leq][y^0] F(x,y,t)$ and $[x^0][y^\leq] F(x,y,t)$, respectively. In these models the negative and non-positive part, respectively, of each axis forms a semipermeable barrier for the walks. In both of these models walks start in the south-east quadrant, they may leave it, but once left, they do not enter it again. Only in the second model walks that end on the origin can neither be extended by a west step nor a south step. Let $F_1=[x^<y^<]F$ and $F_2= F - F_1$. 
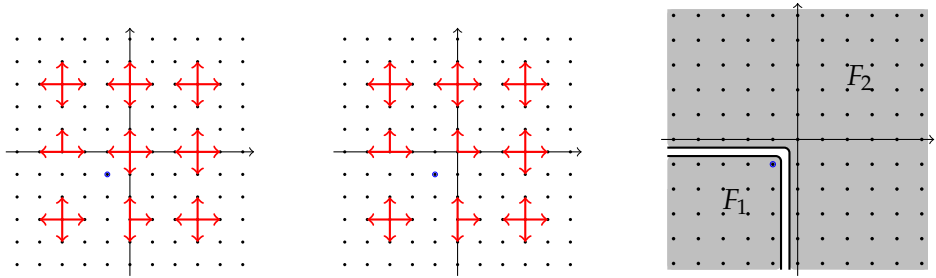
\begin{center}
\begin{tikzpicture}[scale=.3]
  \draw[->](-5.5,0)--(5.5,0);
  \draw[->](0,-5.5)--(0,5.5);
  \foreach \x in {-5,...,5}
  \foreach \y in {-5,...,5}
  \fill (\x,\y) circle(2pt);
  \draw[blue] (-1,-1) circle(3pt);
  \begin{scope}[red,thick]
  \foreach \x in {-3,3}
  \foreach \y in {-3,3}
  \foreach \u/\v in {1/0,0/1,-1/0,0/-1}
  \draw[->] (\x,\y) --+(\u,\v);
  \foreach \x in {-3,3}
  \foreach \u/\v in {1/0,0/1,-1/0} {
    \draw[->] (\x,0) --+(\u,\v);
    \draw[->] (0,\x) --+(\v,\u);
  }
  \draw[->] (3,0)--+(0,-1);
  \draw[->] (0,3)--+(-1,0);
  \draw[->] (0,0)--+(0,1);
  \draw[->] (0,0)--+(0,-1);
  \draw[->] (0,0)--+(1,0);
  \draw[->] (0,0)--+(-1,0);
  \end{scope}
\end{tikzpicture}\hfil
\begin{tikzpicture}[scale=.3]
  \draw[->](-5.5,0)--(5.5,0);
  \draw[->](0,-5.5)--(0,5.5);
  \foreach \x in {-5,...,5}
  \foreach \y in {-5,...,5}
  \fill (\x,\y) circle(2pt);
  \draw[blue] (-1,-1) circle(3pt);
  \begin{scope}[red,thick]
  \foreach \x in {-3,3}
  \foreach \y in {-3,3}
  \foreach \u/\v in {1/0,0/1,-1/0,0/-1}
  \draw[->] (\x,\y) --+(\u,\v);
  \foreach \x in {-3,3}
  \foreach \u/\v in {1/0,0/1,-1/0} {
    \draw[->] (\x,0) --+(\u,\v);
    \draw[->] (0,\x) --+(\v,\u);
  }
  \draw[->] (3,0)--+(0,-1);
  \draw[->] (0,3)--+(-1,0);
  \draw[->] (0,0)--(1,0);
  \draw[->] (0,0)--(0,1);
  \end{scope}
\end{tikzpicture}\hfil
  \begin{tikzpicture}[scale=.33]
    \fill[lightgray,rounded corners=2pt] (-.33,-.33) -- (-5.25,-.33) -- (-5.25,5.25) -- (5.25,5.25) -- (5.25,-5.25) -- (-.33,-5.25) -- cycle
    (-.67,-.67) rectangle (-5.25,-5.25);
    \fill[lightgray] (-2,-2) rectangle (-5.25,-5.25) (-5.25,2) rectangle (5.25,5.25) (2,2) rectangle (5.25,-5.25);
    \draw[->](-5.5,0)--(5.5,0);
    \draw[->](0,-5.5)--(0,5.5);
    \foreach \x in {-5,...,5} \foreach \y in {-5,...,5} \fill (\x,\y) circle(2pt);
    \draw[blue] (-1,-1) circle(3pt);
    \draw[thick,rounded corners=2pt]
    (-5.25,-.67)--(-.67,-.67)--(-.67,-5.25) (-2.5,-2.5) node {$F_1$}
    (-5.25,-.33)--(-.33,-.33)--(-.33,-5.25) (2.5,2.5) node {$F_2$};
  \end{tikzpicture}
\end{center}
As in Section~\ref{sec:all}, we have $F_1(x,y,t) = \bar{x}\bar{y} Q(\bar{x},\bar{y})$ for the $Q(x,y,t)$ from equation~\eqref{eq:Q}. If $F(x,0,t) = [x^<y^0]F(x,y,t)$ and $F(0,y,t) = [x^0y^<]F(x,y,t)$ (the
left-most case on the previous figure) we can show D-finiteness of $F_2$ using the analytic method~\cite{raschel12, raschel19}.
The argument presented in \cite{raschel19} is based on the decomposition $F_2=F_2^U+F_2^D+F_2^L$ with
$F_{2}^{D}=\sum_{i \geq 0} x^{i}y^{i} [x^{i} y^{i}]F_{2}$ and $F_{2}^{L}=\sum_{i \geq 0, j \leq i-1}x^{i}y^{j}[x^{i}y^{j}]F_{2}$ and $F_{2}^{U}(x,y,t) = F_{2}^{L}(y,x,t)$.
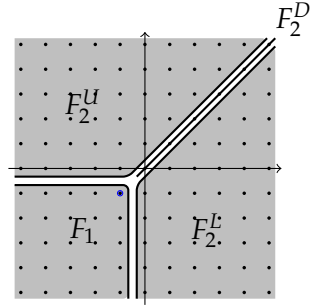
\begin{center}
  \begin{tikzpicture}[scale=.33]
    \fill[lightgray,rounded corners=2pt] 
    (-.33,-5.25)--(-.33,-.67)--(5.25,4.92)--(5.25,-5.25)
    (-5.25,-.33)--(-.67,-.33)--(4.92,5.25)--(-5.25,5.25)
    (-5.25,-.67)--(-.67,-.67)--(-.67,-5.25)--(-5.25,-5.25);
    \draw[->](-5.5,0)--(5.5,0);
    \draw[->](0,-5.5)--(0,5.5);
    \foreach \x in {-5,...,5} \foreach \y in {-5,...,5} \fill (\x,\y) circle(2pt);
    \draw[blue] (-1,-1) circle(3pt);
    \draw[thick,rounded corners=2pt]
    (-5.25,-.67)--(-.67,-.67)--(-.67,-5.25) (-2.5,-2.5) node {$F_1$}
    (-.33,-5.25)--(-.33,-0.67)--(5.25,4.92) (2.5,-2.5) node {$F_{2}^{L}$}
    (-5.25,-.33)--(-.67,-.33)--(4.92,5.25) (-2.5,2.5) node {$F_{2}^{U}$}
    (-.33,-.33)--(5.25,5.25) (6,6) node {$F_{2}^{D}$};
  \end{tikzpicture}
\end{center}
The functions $F_2^D$ and $F_2^L$ satisfy the equations
\begin{alignat*}5
F_{2}^{D} &= 2 t ( \bar x + y ) \sum_{i \geq 0} x^{i} y^{i-1} [x^{i} y^{i-1}] F_{2}^{L}  &&  && - 2 t\bar x \bar y [x^{0} \bar y] F_{2}^{L}\\
F_{2}^{L} &= \smash{\underbrace{\rule[-5pt]{0pt}{0pt} t [\bar x] F_{1} + t ( x + \bar y ) F_{2}^{D}\qquad\quad\qquad}_{\text{\mathstrut initial conditions}}}&& + \smash{\underbrace{S t F_{2}^{L}}_{\text{\mathstrut recurrence}}} && 
- t ( \bar x + y ) \sum_{i \geq 0} x^{i} y^{i-1}[x^{i} y^{i-1}] F_{2}^{L}  \\
  & && &&\underbrace{\rule[-5pt]{0pt}{0pt} + t \bar x \bar y[x^{0} \bar y] F_{2}^{L}  - t \bar x [x^{0}] F_{2}^{L} \quad\quad} _{\text{\mathstrut boundary conditions}}.
\end{alignat*}
Eliminating the term involving the infinite sum gives the equation
\begin{equation}\label{star}\tag{$\ast$}
 (1-S t) F_2^L = t [\bar x] F_1 + (t x + t\bar y - \tfrac{1}{2}) F_2^D - t \bar x[x^{0}] F_{2}^L.
\end{equation}
From here on, we can follow the derivation in~\cite{raschel19} step by step and obtain an expression of
the form 
\[
F_2^D(xy,\bar x, t)=A(y)\oint B(y,z)[\bar y]F_1\bigl(yC(z),1/C(z), t\bigr)dz,
\]
where $A$, $B$, and $C$ are certain algebraic functions, the integral is taken around the unit circle, and $t$
is viewed as a fixed small positive real number.
The interested reader will find on our website~\cite{buchacher20a} a Maple session in which this derivation
is worked out in full detail, and where also explicit expressions for $A$, $B$, and $C$ are provided.
What matters here is that the D-finiteness of $F_1$ together with the algebraicity of $A,B,C$ implies the
D-finiteness of~$F_2^D$. Knowing this we can solve the functional equation \eqref{star} for $[x^0]F_2^L$
after setting $x$ to a root $X(y,t)\in\set Q[y,\bar y][[t]]$ of the polynomial $1-St\in\set Q[x,\bar x,y,\bar y,t]$
in order to eliminate the left hand side. The resulting expression
\[
  [x^0]F_2^L = X(y,t) [\bar x]F_1 + X(y,t)(X(y,t) + \bar y - \frac{1}{2t}) F_2^D(X(y,t),y, t)
\]
certifies that $[x^0]F_2^L$ is D-finite. With the knowledge that $F_1$, $F_2^D$, and $[x^0]F_2^L$ are D-finite,
it follows from \eqref{star} that $F_2^L$ is D-finite. Then $F_2^U(x,y,t)=F_2^L(y,x,t)$ is D-finite as well,
and it finally follows that $F_2=F_2^D+F_2^L+F_2^U$ and $F=F_1+F_2$ are D-finite, as claimed.

Like in the previous section, we were not able to construct a certified annihilating operator for the series $F(1,1,t)$
but only have a convincing guess. Assuming however that this guess is correct, we can again show that the series $F(1,1,t)$ must
be transcendental. The reasoning is like in the previous section: the guessed operator now has order 10 and can be written
as the least common left multiple of four irreducible operators, exactly one of them admits transcendental
solutions and therefore only has transcendental solutions. The lclm of the remaining operators does not annihilate $F(1,1,t)$, so $F(1,1,t)$ must be transcendental.
The operators are available on our website~\cite{buchacher20a}.

\bigskip
Unfortunately, we are not able to prove D-finiteness of the solution for the interpretation of $F(0,y,t)$
and $F(x,0,t)$ as $[x^0][y^\leq] F(x,y,t)$ and $[y^0][x^\leq]F(x,y,t)$, respectively. If we proceed as above,
we are led to an expression of the form 
\[
F_2^D(xy,\bar x, t)=A(y)\oint B(y,z)\Bigl([\bar y]F_1\bigl(yC(z),1/C(z), t\bigr) - [y^0]F_2^D(xy,\bar x,t)\Bigr)dz,
\]
where $A$, $B$, and $C$ are again algebraic functions, the integral is taken around the unit circle, and $t$
is viewed as a fixed small positive real number. As we do not know whether $[y^0]F_2^D(xy,\bar x,t)$ is D-finite,
we are stuck at this point.
It does seem however that the solution $F(x,y,t)$ is D-finite also in this case, at least for $x=y=1$.
We have found a convincing candidate for an annihilating operator of order 13 by guessing.
Again, this operator can be written as the least common left multiple of irreducible operators $L_1,\dots,L_6$,
available on our website, so we can write $F(1,1,t)=f_1 + \cdots + f_6$ for certain solutions $f_i$ of the irreducible
right factors~$L_i$. The difference to the earlier cases is that now several of the $L_i$ have transcendental solutions,
so in order to show that $F(1,1,t)$ is not algebraic, we must show that their sum is not algebraic.
This can be done by constructing an operator $P$ which annihilates all but one of the
transcendental summands~$f_i$, so that $P\cdot F(1,1,t)$ can be written as a sum of one transcendental series and some
algebraic series, which asserts that $P\cdot F(1,1,t)$ is transcendental. But then $F(1,1,t)$ must be transcendental as well,
because if it were algebraic, so would be $P\cdot F(1,1,t)$. 

\bigskip
In summary, we have shown that for the interpretation $F(0,y,t)=[x^0][y^<] F(x,y,t)$ and $F(x,0,t)=[y^0][x^<] F(x,y,t)$,
the solution $F(x,y,t)$ of the functional equation~\eqref{eq:F} is D-finite. Moreover, under the hypothesis that guessed
annihilating operators for the solutions at $x=y=1$ are correct, we can also show that these series are transcendental.
For the interpretation $F(0,y,t)=[x^0][y^\leq] F(x,y,t)$ and $F(x,0,t)=[y^0][x^\leq] F(x,y,t)$, we have no proof that the
solution $F(x,y,t)$ is D-finite, but we have a guessed equation for $F(1,1,t)$ whose correctness implies the transcendence
of the solution.

\section{No Compartments}\label{sec:pp}

One other variant of \eqref{eq:F} is when we read $F(0,y,t)$ as $[x^0][y^>] F(x,y,t)$, and $F(x,0,t)$ likewise.
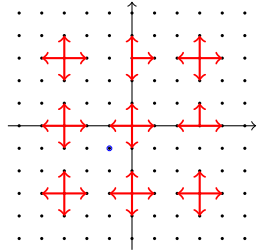
\begin{center}
\begin{tikzpicture}[scale=.3]
  \draw[->](-5.5,0)--(5.5,0);
  \draw[->](0,-5.5)--(0,5.5);
  \foreach \x in {-5,...,5}
  \foreach \y in {-5,...,5}
  \fill (\x,\y) circle(2pt);
  \draw[blue] (-1,-1) circle(3pt);
  \begin{scope}[red,thick]
  \foreach \x in {-3,3}
  \foreach \y in {-3,3}
  \foreach \u/\v in {1/0,0/1,-1/0,0/-1}
  \draw[->] (\x,\y) --+(\u,\v);
  \foreach \x in {-3,3}
  \foreach \u/\v in {1/0,0/1,-1/0} {
    \draw[->] (\x,0) --+(\u,\v);
    \draw[->] (0,\x) --+(\v,\u);
  }
  \draw[->] (-3,0)--+(0,-1);
  \draw[->] (0,-3)--+(-1,0);
  \draw[->] (0,0)--+(0,1);
  \draw[->] (0,0)--+(0,-1);
  \draw[->] (0,0)--+(1,0);
  \draw[->] (0,0)--+(-1,0);
  \end{scope}
\end{tikzpicture}
\end{center}
For this model there is no natural division of the domain into compartments that cannot be left once entered.
Even though a walk may pass through the positive part of either axis only in one direction, it can still escape from
the first quadrant through the origin.
It is arguably for this reason that this model appears to be more difficult than the others.
Indeed, we have not been able to solve it.

Computer experiments with the first 2000 series terms suggest that the coefficient sequence of $F(1,1,t)$ grows asymptotically
like $c 4^n n^{-1/3}$ for $n\to\infty$ and some constant $c\approx 1.91$. Moreover, for the number $a_n$ of walks of length $2n$
starting and ending at $(-1,-1)$, we find, based on 6300 sequence terms, a conjectured asymptotic behaviour of the
form $c 4^n n^{-5/3}$ for $n\to\infty$ and some nonzero constant~$c$.
Even if these growth rates are correct, they cannot even be used to exclude algebraicity of the generating functions.

We have also searched for candidates for algebraic and differential equations by guessing based on almost 98000 sequence terms
(modulo~45007), but did not find any. This implies that such equations, if they exist, must be extraordinarily large. We are
tempted to conjecture that they do not exist, i.e., that the solution $F(x,y,t)$ of \eqref{eq:F} is not D-finite for the interpretation
under consideration. The terms we computed can be found on our website. 

\section{Conclusion}

We investigated the functional equation
\begin{equation*}
(1 - (x+y+\bar x+\bar y)t) F(x,y,t) = \bar x\bar y - \bar x t F(0,y,t) - \bar y t F(x,0,t)
\end{equation*}
and its solution $F(x,y,t)$ in $\mathbb{Q}[x,\bar{x},y,\bar{y}][[t]]$ for different interpretations of $F(0,y,t)$ and $F(x,0,t)$. For $F(0,y,t) = [x^0]F$ and $F(x,0,t)=[y^0]F$ we answered the main questions: we proved that $F$ is D-finite, and we showed that it is not algebraic. As we recall in the introduction, there is an equivalence between zeroness of the orbit sum and algebraicity of the generating function. However, this result can only be applied if the generating functions can be meromorphically continued on the universal covering~\cite[Sec.~3]{kurkova15}, which is not the case for $F(x,y,t)$ nor any interpretation of $F(x,0,t)$ and $F(0,y,t)$. So our result is not in conflict with these earlier results. 

For other interpretations there are several open points which would deserve further consideration. One point is the pending proof of the guessed operators on which the transcendence arguments of Sections~\ref{sec:nnp}--\ref{sec:np/npp} rely.
In the second case considered in Section~\ref{sec:np/npp}, not only the guessed operator but also D-finiteness in general remains to be proven.
Another open issue is the clarification of the nature of the solution in Section~\ref{sec:pp}: is it really non-D-finite?
Besides answering these open questions, there are some natural extensions and generalizations which could be addressed.
For example, we have only considered analogous interpretations for $F(x,0,t)$ and $F(0,y,t)$ in this paper, but mixed cases such as $F(x,0,t)=[y^0][x^>]F(x,y,t)$
and $F(0,y,t)=[x^0][y^\leq]F(x,y,t)$ might also be interesting. First experiments suggest that some combinations are D-finite.
Another possible variation concerns the start point. There are other points besides $(-1,-1)$ which lead to a zero orbit sum, for
example $(-1,1)$. Can the starting point affect the nature of the solution?
Finally, we have restricted ourselves to the case of simple walks, and it would be interesting to see what happens for other step sets.

\printbibliography
\end{document}